\documentclass[a4paper,UKenglish,cleveref, autoref, thm-restate]{lipics-v2021}
\graphicspath{{./graphics/}}

\usepackage[svgnames]{xcolor}

\title{Waterproof Editor: an educational environment for proof assistants and programming languages} 

\titlerunning{Waterproof Editor} 

\author{Pim Otte}{Mathematical Institute, Utrecht University, The Netherlands \and Department of Mathematics and Computer Science, Eindhoven University of Technology, The Netherlands}{p.j.otte@uu.nl}{https://orcid.org/0009-0001-4076-0616}{}

\author{Dick Arends}{Radboud University Nijmegen, The Netherlands}{dick.arends@ru.nl}{https://orcid.org/0009-0009-3685-1572}{}

\author{Raul S\'anchez Flores}{Department of Mathematics and Computer Science, Eindhoven University of Technology, The Netherlands}{r.sanchez.flores@student.tue.nl}{}{}

\author{Pieter Wils}{Department of Mathematics and Computer Science, Eindhoven University of Technology, The Netherlands}{p.f.e.j.wils@student.tue.nl}{}{}

\author{Jim Portegies}{Department of Mathematics and Computer Science, Eindhoven University of Technology, The Netherlands}{j.w.portegies@tue.nl}{https://orcid.org/0000-0002-2103-7334}{}

\authorrunning{Otte, Arends, S\'anchez Flores, Wils, Portegies} 

\Copyright{Pim Otte and Dick Arends and Raul S\'anchez Flores and Pieter Wils and Jim Portegies} 



\category{} 

\relatedversion{} 

\acknowledgements{We would like to thank Thom Coolen, Piotr Czarnecki, Jorge Dias Batista, Stefan Dragoslavov, Oliwia Hejduk, Serhiy Medvedyev, Vincente de Oliveira Vegar, Sinjini Pande, Militsa Stefanova and Marcel \.{Z}elazny for their work in their Software Engineering Project that allowed the use of Lean for Waterproof as a backend. We would like to thank Mihail Abramov, Benjamin Blankers, Tammo Steffens and Lars Wagoner for their contributions in the Waterproof team. We thank Jelle Wemmenhove for their contributions to Waterproof and useful suggestions on educational matters. We are very grateful to Emilio J. Gallego Arias and Shachar Itzhaky for their support and in particular for realizing a webworker for Waterproof. We would like to thank 4TU.AMI+ and the Drive program at Eindhoven University of Technology for financial support.}

\EventEditors{}
\EventNoEds{0}
\EventLongTitle{Tools for Educational Activities in Logic}
\EventShortTitle{TEAL 2026}
\EventAcronym{TEAL}
\EventYear{2026}
\EventDate{July 26, 2026}
\EventLocation{Lisbon Portugal}
\EventLogo{}
\SeriesVolume{}
\ArticleNo{20}

\hideLIPIcs
\nolinenumbers
\begin{document}

\maketitle

\begin{abstract}
Waterproof Editor provides an educational environment specifically targeted to teaching with proof assistants or programming languages. It arose from Waterproof, educational software targeted at helping students acquire the skill of giving mathematical proofs. Its original features such as enabling rich formatting and providing clear input areas are now abstracted away in an npm package and can be used in different educational contexts. 
We invite interested parties to use this component in their educational software, and offer to assist with this.
\end{abstract}

\section{Introduction}

We propose to demo \emph{Waterproof Editor}, an environment that was originally designed to work with proof assistants in education.
Waterproof Editor was initially a tightly interwoven part of Waterproof \cite{wemmenhove2024waterproof}, an educational tool aimed at helping students learn to write mathematical proofs. Waterproof was built on top of the proof assistant Rocq, and leveraged this to provide students with immediate feedback.

\begin{figure}
    \centering
    \includegraphics[width=0.5\linewidth]{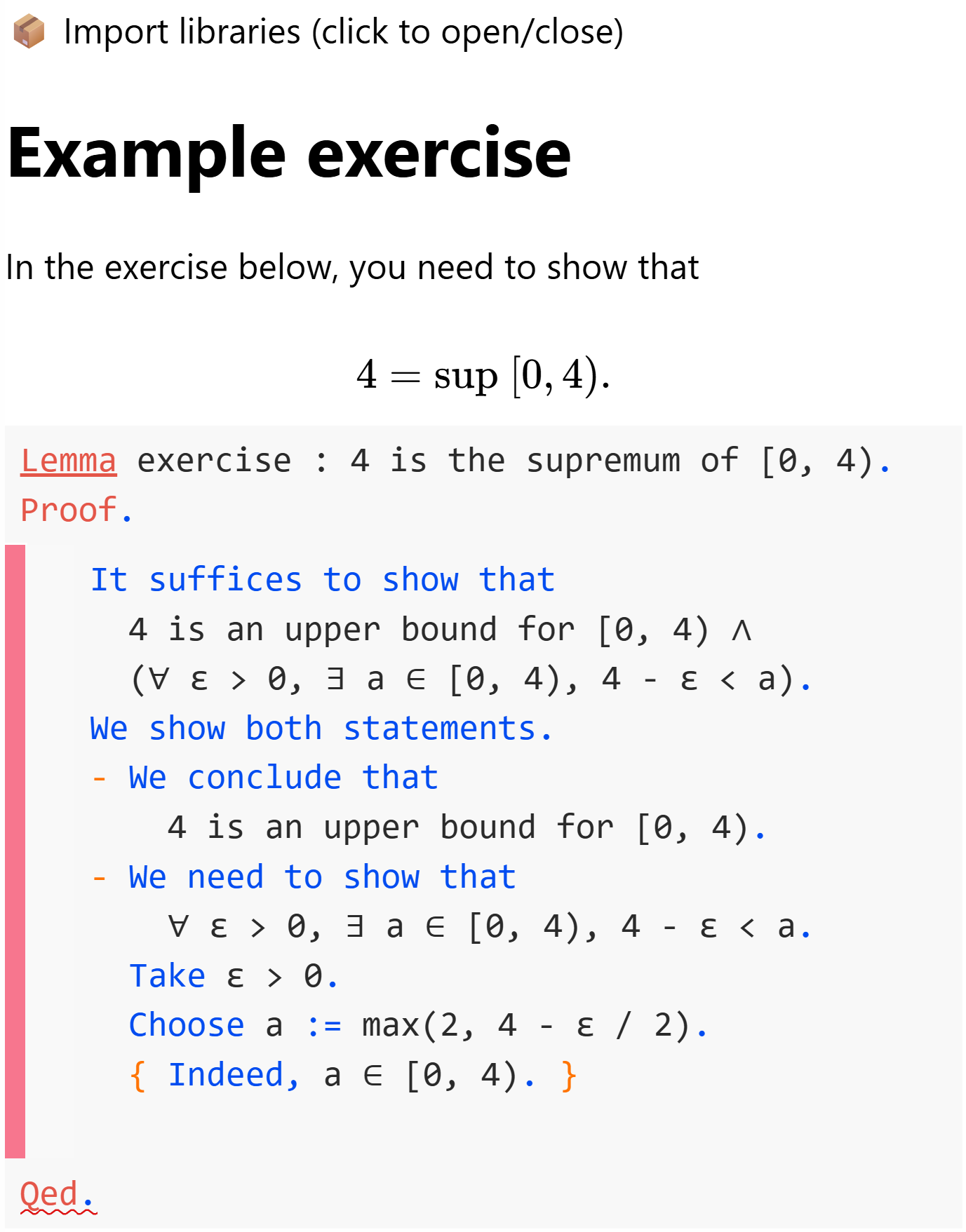}
    \caption{A screenshot of Waterproof with an input area designated by the red bar on the left. The bar is red as the proof is not yet complete.}
    \label{fi:screenshot-waterproof}
\end{figure}

In Waterproof, students can write their proofs in designated input areas, the input areas are colored on the basis of evaluation of the proofs, accidental edits outside of these input areas are prevented, diagnostic messages are clearly presented underneath input blocks, and Markdown is rendered to provide documents with a structured look, see Figure \ref{fi:screenshot-waterproof} for a screenshot.
From a desire to also have such an environment available for working with other proof assistants, in particular Lean, we split off large parts of this functionality into the Waterproof Editor project.

With abstracting a lot of the features in Waterproof to Waterproof Editor, the environment suddenly also becomes useful for many other purposes, for instance for working with programming languages.
Conceptually, the only requirement is an evaluation function on the content of the input areas. 
In this demo we will present several use cases, in which these evaluation functions determine whether a proof is correct or whether JavaScript code passes a suite of tests.

Furthermore, since Waterproof Editor is distributed as an npm package, it helped us realize a long-standing wish of creating a browser-only version of Waterproof, in which teachers or students can try out or even work with Waterproof without any installation by just following a link.

We see it as a large advantage that Waterproof Editor allows teachers to provide assignments to students in which the solutions can be put in the same document as the exercise statements.

In the demo, we will show several use cases, including a version of Waterproof for Lean, programming exercises in JavaScript and the browser version of Waterproof.
We extend an open invitation for people to imagine and implement new applications, for which we would gladly provide the help.

\subsection{Related work}

With the respect of mixed prose-code documents, several projects are closely related.
Alectryon~\cite{pitclaudel2020untangling} combines the mixed document format with
full elaboration of either Lean or Rocq code, yielding stand-alone documents with elaborated proof states. jsCoq~\cite{arias2017jscoq} provides an interactive standalone browser environment,
and pioneered the web-worker infrastructure we reuse for our use case of Waterproof in the browser (see Section \ref{sec:browser}). Jupyter Notebooks~\cite{jupyter} share conceptual and visual similarities, but
generally allows out-of-order execution, which is unsuitable for our use case. ProseMirror\footnote{\url{https://prosemirror.net/}} and CodeMirror\footnote{\url{https://codemirror.net/}} are the major components used to realize Waterproof Editor.

Diproche~\cite{carl2020number} provides a controlled natural language environment built on top of a custom automated theorem prover, while Verbose Lean~\cite{massot2024teaching} uses Lean. For the use case ``Waterproof with Lean'' (Section \ref{sec:lean}) we have combined Verbose Lean with Verso\footnote{\url{https://verso.lean-lang.org/}}, a documentation system for Lean that allows for documents to be written alongside the formalization.

\section{Background on Waterproof}

Waterproof has been developed at Eindhoven University of Technology~\cite{wemmenhove2024waterproof} and used in the Analysis 1 course 
since 2019.
In 2025, Waterproof has been piloted in the introductory proof course at Utrecht University.
For a more complete overview of (the history) of Waterproof we refer to~\cite{portegies2025waterproof}.

Waterproof helps students learn the skill of writing mathematical proofs.
Waterproof presents the feedback from a proof assistant (Rocq or Lean) and
uses a layer of controlled natural language (coq-waterproof\footnote{\url{https://github.com/impermeable/coq-waterproof}} or Verbose Lean \cite{massot2024teaching})
to ease the transfer of skills obtained to writing paper proofs.

Previous work has evaluated the use of Waterproof in education through the lens of service design
\cite{wemmenhove2025service}, as well as the effect on the structure of the proofs that students deliver \cite{wemmenhove2026comparative}.
The latter showed that the proof structure used by Waterproof users differed in three ways from the structure used by non-Waterproof users.
Waterproof users did introduce variables more often and they announced the goal to be proven more frequently, but at the same time they were more likely to forget to check domain conditions in quantifiers.
This third finding has been traced to a specific aspect of Waterproof, and this has since been improved.

\section{Features of Waterproof Editor}

All features of Waterproof Editor have been designed for the purpose of (mathematics) education. 
Markdown rendering provides students with nicely formatted exercise sheets, increasing readability. 
Collapsible areas have a dual use for static hints, which students can reveal on demand, and technical details, which might be unfolded by
a teacher for debugging purposes.
Input areas prevent the students from accidentally editing the rest of the document, which includes exercise statements.
Feedback and messages are displayed in the document to decrease the amount of different areas students need to pay attention to.
Visual editing of the document allows teachers to interactively check the exercise sheets as if they were a student,
and keeping the option of plaintext editing ensures that teachers who are power-users of the underlying proof assistant can
edit exercise sheets with their usual workflow.

\subsection{Markdown rendering} 

To give the documents an attractive and structured look, we make use of Markdown rendering. The screenshot in Figure \ref{fi:screenshot-waterproof} shows a bold face header, accompanying text and a rendered \LaTeX{} formula. It is also possible to embed pictures and videos.
We believe that this rendering has many advantages instead of presenting text-based proof documents. Certainly, the resulting rendered document has many similarities with Jupyter Notebooks, while the underlying document is very different.

\subsection{Collapsible areas for hints and technical details}

Waterproof Editor allows for collapsing certain areas of the document so that they are not visible to the user by default. This is for instance useful when providing an optional hint to the students, that only unfolds when they click it, or when one would not want to distract the user with technical details. In the document in Figure \ref{fi:screenshot-waterproof}, it is used to hide all the statements that import the libraries and set up the automation for the exercise sheet to work.

\subsection{Input areas}

Students write their solutions in designated input areas. In Figure \ref{fi:screenshot-waterproof}, the red bar next to the input area indicates in this case that the proof is not complete yet. In this particular example of Waterproof, the bar would also color red if the solution until that point would contain an error.
If the input block would only contain a warning, the bar would be colored orange.
If the evaluation is not complete yet, the bar would be blue.
The exact meaning of the colors in Waterproof Editor is customizable: one only needs an evaluation function.
We note that this evaluation function can be reused in (server-side) automatic grading of exercises.

A main advantage of specifying areas in which students can edit is that we prevent them from accidentally changing anything in the exercise sheet outside of input areas.
Our experience in another large course at Eindhoven University of Technology relying on Jupyter Notebooks for Python programming was that if one does not prevent such accidental edits, they will invariably happen, rendering exercise sheets useless either for the students themselves or later when they have to submit them.

\subsection{Feedback and messages below codeblocks}

\begin{figure}
\centering
    \includegraphics[width=1\linewidth]{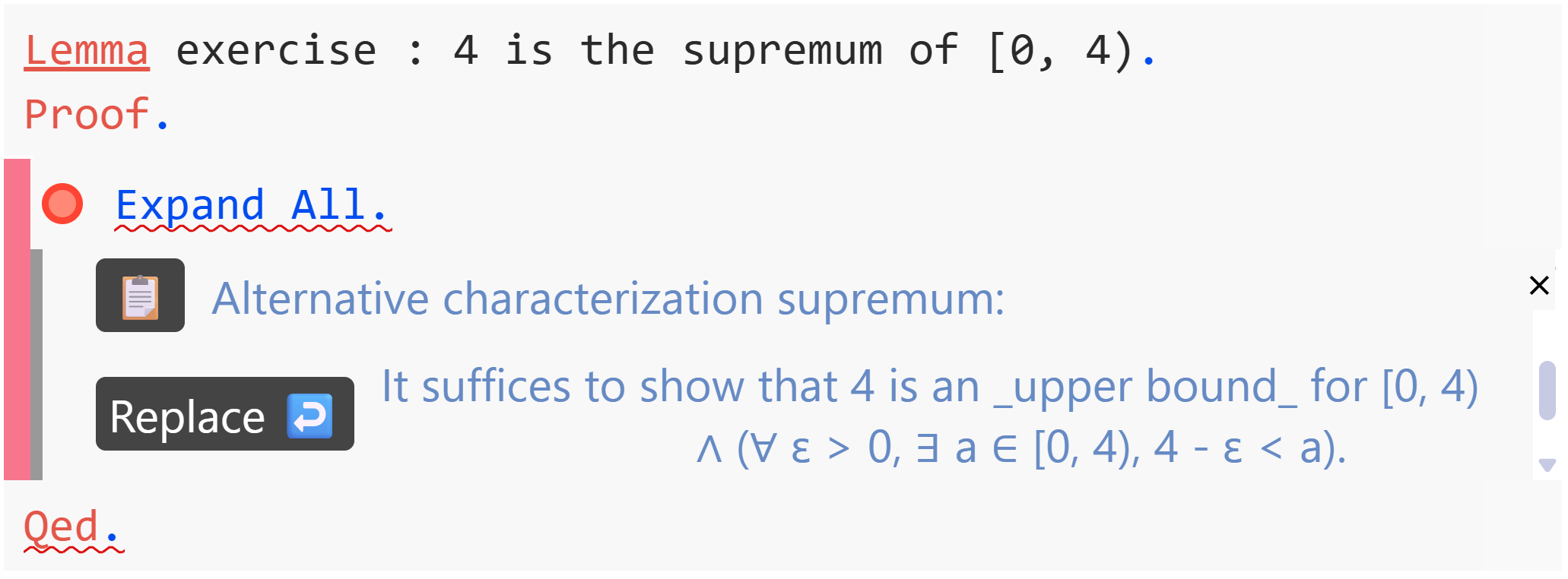}
    \caption{A codeblock with suggestions on what to do next. Pressing the ``Replace'' button would replace the line with \texttt{Expand All.} by the suggested alternative characterization.}
    \label{fi:screenshot-feedback}
\end{figure}

We rely on a feature of CodeMirror that allows for putting diagnostic messages directly underneath codeblocks.
We use this feature to communicate error messages and warnings, but also to give students suggestions on what to do next, including buttons that allow students to carry through fixes or suggestions with a single click, see Figure \ref{fi:screenshot-feedback}.

\subsection{Visual and plaintext editing}

An important core design principle is that Waterproof documents should be editable for students and teachers in both Waterproof and a plaintext format, using the standard tooling.\footnote{An unexpected instance where this proved useful was when a student requested whether they ``could work in Emacs, because VS Code is non-free''} Standard Waterproof documents are \texttt{.mv} files, which are effectively Markdown files in which the code blocks contain Rocq code.
Expert teachers can then easily manipulate an exercise sheet in plaintext mode, while if teachers want to immediately have a feeling on what an exercise sheet would look like if they render it, they can edit in Waterproof.

\section{Use case: Waterproof with Lean}\label{sec:lean}

\begin{figure}
    \centering
    \includegraphics[width=0.8\linewidth]{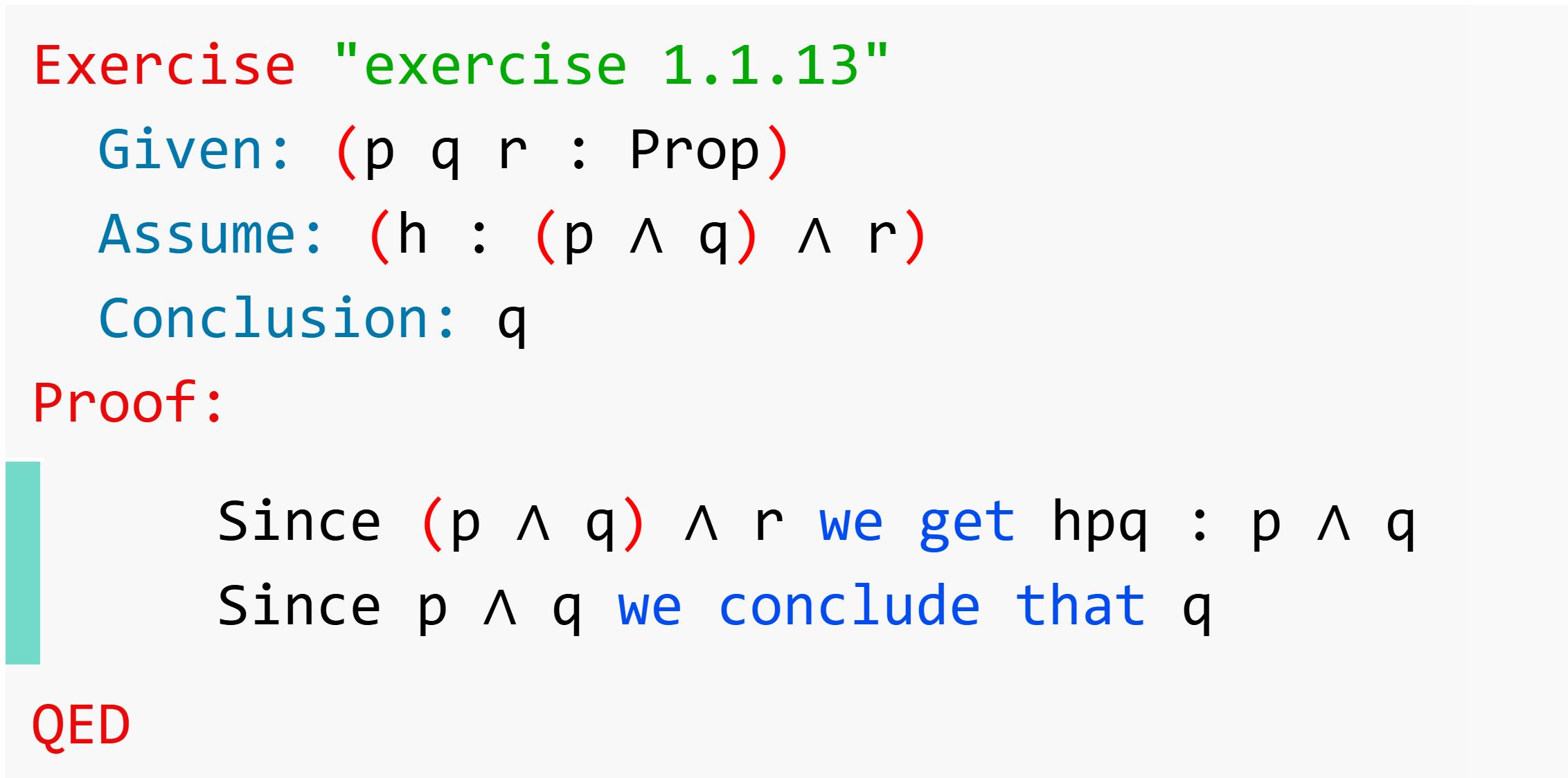}
    \caption{A screenshot of Waterproof Editor with a Verbose lean exercise.}
    \label{fi:waterproof-lean}
\end{figure}

As a first use case for Waterproof Editor, we present support for Verbose Lean~\cite{massot2024teaching} in Waterproof, see Figure \ref{fi:waterproof-lean}. It provides the same editing functionality as for the original Rocq-based exercise sheets.

Of course, to make this work, also changes to the Waterproof VS Code extension were necessary, but splitting of Waterproof-editor did make this significantly easier.
A large part of this extra work was performed as a Software Engineering Project at the Eindhoven University of Technology. The Waterproof goal window was replaced using the `lean4-infoview' package and some of the Rocq LSP client functionality could be abstracted and reused. We implemented a custom Verso Genre to ensure the plain files can be edited with the standard tools; Verso provides full support of Markdown-like syntax, with the exception of a preamble that needs special treatment when parsing.

\section{Use case: teaching programming using JavaScript}

\begin{figure}
    \centering
    \includegraphics[width=\linewidth]{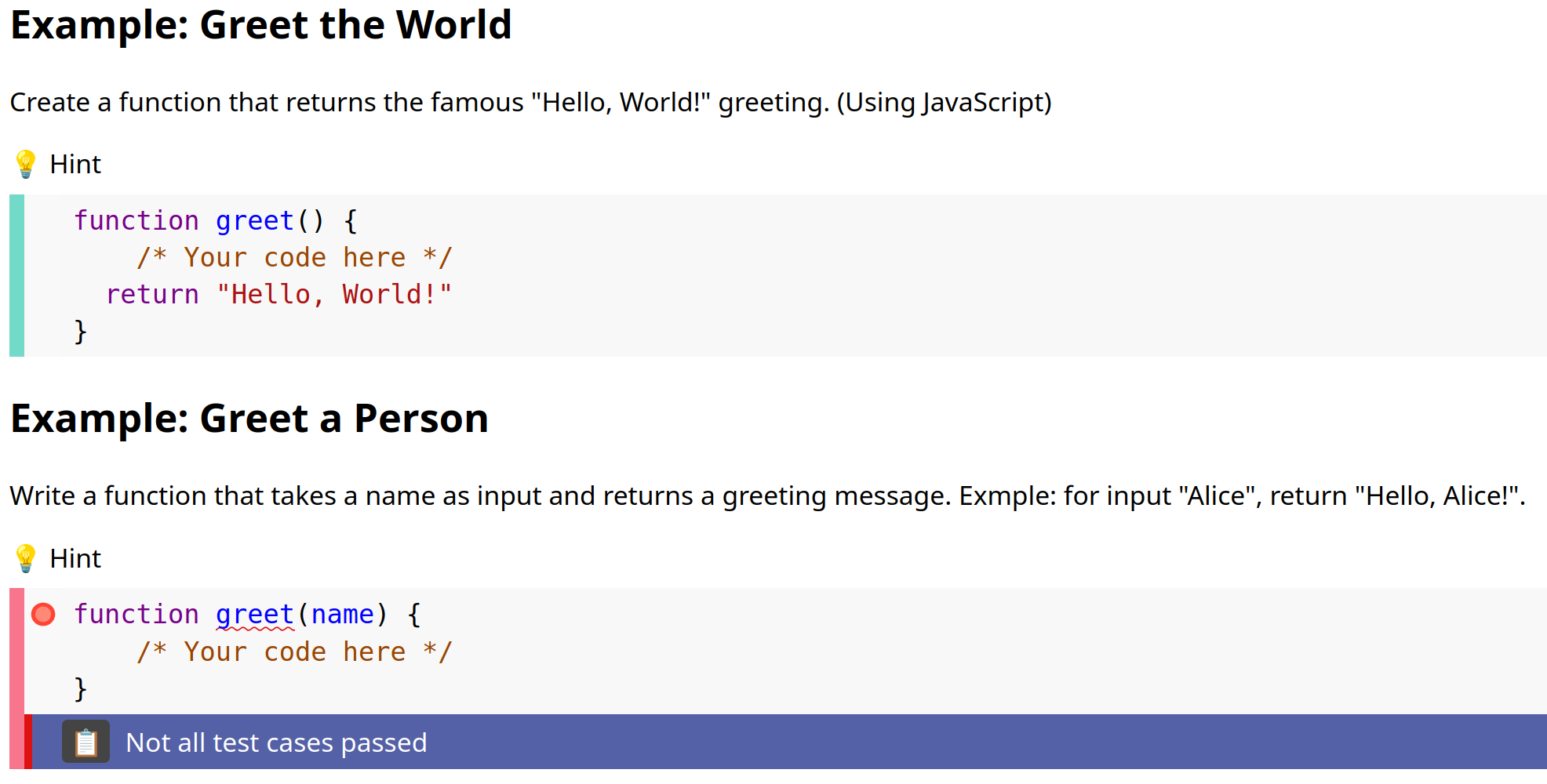}
    \caption{Illustration of JavaScript exercises using Waterproof Editor.}
    \label{fi:screenshot-javascript}
\end{figure}

The use of Waterproof Editor is not limited to proof assistants.
As a validation of the functionality, we implemented a proof of concept\footnote{Live version at \url{https://impermeable.github.io/waterproof-javascript/}, source code at \url{https://github.com/impermeable/waterproof-javascript}} using the editor for a different educational context, namely for teaching algorithms and programming with JavaScript, see Figure \ref{fi:screenshot-javascript}.

In this case, the colors next to input areas get different meanings.
An input area turns green if all test cases have passed, it turns red if one or more test cases fail. If that happens, a diagnostic error message underneath the input block informs the student. The bar next to the input block turns yellow if the execution exceeded a predetermined time limit. In that case, an ``Execution timed out'' warning message appears.

\section{Use case: browser version of Waterproof}\label{sec:browser}

Another advantage of Waterproof Editor is that it is easy to use it in a different environment, such as in a browser environment.
For years, it has been an objective of the Waterproof team to make a version of Waterproof that runs fully in the Browser.
We got much closer to this objective with the release of a WASM webworker for Rocq-lsp, a Rocq LSP server.
Relying on this webworker, in the Summer of 2025 we released a version of Waterproof that was able to run on VSCode.dev. Note that no external server was necessary for this to run, and students did not need to install anything.

The easiest way to work with this requires students to have a GitHub account. Working with a local folder is also possible, and it only takes a few lines of instructions to get started with Waterproof in this way.

By building on Waterproof Editor, we could build a version of Waterproof that runs fully in the browser\footnote{Live version at \url{https://impermeable.github.io/aquarium}, source code at \url{https://github.com/impermeable/aquarium}}, and does not rely on the web version of Visual Studio Code.
This significantly lowers the barrier for students to get started, as well as for prospective educators who would want to test Waterproof.
They have access to several exercise repositories to immediately start working with Waterproof code.

\section{What is needed to use Waterproof Editor for a custom application?}

In order to extend Waterproof with another language, all that is needed to provide a language configuration, as well as a configuration for highlighting, and possibly adapting one of the existing parsers to parse the exact on-disk format at the time the doucment is first opened.
For proof assistants, which generally have some amount of non-standard LSP functionality, some glue code is also required.

Embedding Waterproof Editor in a similar way as the standalone version in the browser requires providing a configuration, which contains a set of callbacks that allow
the editor to communicate to the application.

Either of these options are likely doable by student projects. Implementation of the presented use cases required finding the rough edges for integration and allows future projects to integrate only by implementing specific details. Tooling that already supports a literate Markdown-like format is likely very easy to integrate by reusing even more of the available infrastructure. 

\section{Conclusion}

We believe that the breadth of the presented use cases shows the wide applicability of Waterproof Editor as a component in teaching tools. We hope that we can demo the tool at TEAL and excite some people to try to use Waterproof Editor for their own application.

\bibliography{main}

\end{document}